\documentclass[a4paper]{ifacconf}
\usepackage{natbib}            % you should have natbib.sty
\usepackage{graphicx}          % Include this line if your
\usepackage[dvips]{epsfig}     % or this line, depending on which
\usepackage{dahlin,enumitem}

\usepackage[latin1]{inputenc}
\usepackage{amsmath}
\usepackage{amssymb}
\usepackage{epsfig}
\usepackage{psfig}
\usepackage{psfrag}
\usepackage{graphics}
\usepackage{color}
\usepackage{pmat}
\usepackage{subfigure}
\newcommand{\Ebf}{{\mathbf E}}

\newcommand{\xbf}{{\mathbf x}}

\newcommand{\zbf}{{\mathbf z}}

\newcommand{\Ccal}{{\mathcal {C}}}

\newcommand{\Gcal}{{\mathcal {G}}}
\newcommand{\Ical}{{\mathcal {I}}}

\newcommand{\Pcal}{{\mathcal {P}}}

\newcommand{\Vcal}{{\mathcal {V}}}
\newcommand{\vcal}{{\mathcal {V}}}

\newcommand{\rb}{{\mathbb R}}

\newcommand{\beq}{\begin{equation}}
\newcommand{\eeq}{\end{equation}}

\newcommand{\bnul}{\begin{enumerate}[a)]}
\newcommand{\enul}{\end{enumerate}}

\newcommand{\op}[1]{\operatorname{#1}}

\newcommand{\fin}{\hspace*{\fill}~$\blacksquare$}

\allowdisplaybreaks

\begin{document}

\begin{frontmatter}

\title{A graph/particle-based method for experiment design in nonlinear systems\thanksref{footnoteinfo}} % Title, preferably not more than 10 words.

\thanks[footnoteinfo]{This work was supported in part by the Swedish Research Council under contracts 621-2011-5890 and 621-2009-4017, and in part by the European Research Council under the advanced grant LEARN, contract 267381.}

\author[First]{Patricio E. Valenzuela}
\author[Second]{Johan Dahlin}
\author[First]{Cristian R. Rojas}
\author[Third]{Thomas B. Schön}

\address[First]{Department of Automatic Control, KTH Royal Institute of Technology, SE-100 44, Stockholm, Sweden (e-mail: \{pva, crro\}@kth.se).}
\address[Second]{Division of Automatic Control, Linköping University, SE-581 83 Linköping, Sweden (e-mail: johan.dahlin@liu.se).}
\address[Third]{Division of Systems and Control, Uppsala University, SE-751 05 Uppsala, Sweden (e-mail: thomas.schon@it.uu.se).}

\begin{keyword}                           % Five to ten keywords,
System identification, input design, particle filter, nonlinear systems.               % chosen from the IFAC
\end{keyword}                             % keyword list or with the
                                          % help of the Automatica
                                          % keyword wizard

\begin{abstract}                          % Abstract of not more than 250 words.
We propose an extended method for experiment design in nonlinear state space models. The proposed input design technique optimizes a scalar cost function of the information matrix, by computing the optimal stationary probability mass function (pmf) from which an input sequence is sampled. %, extending the class of input sequences considered in previous results for state space models.
 %Since
 The feasible set of the stationary pmf is a polytope, allowing it to be expressed as a convex combination of its extreme points. %By assuming a finite set of possible values for the input, we can compute the extreme points in the feasible set using graph theory.
 The extreme points in the feasible set of pmf's can be computed using graph theory. Therefore, the final information matrix can be approximated as a convex combination of the information matrices associated with each extreme point. For nonlinear systems, the information matrices for each extreme point can be computed by using particle methods. %, which extends the set of feasible model structures covered by an existing input design method.
 Numerical examples show that the proposed technique can be successfully employed for experiment design in nonlinear systems.
\end{abstract}

\end{frontmatter}

%%%%%%%%%%%%%%%%%%%%%%%%%%%%%%%%%%%%%%%%%%%%%%%%%%%%%%%%%%%%%%%%%%%%%%%%%%%%%%%%
\section{Introduction}
Experiment design deals with the generation of an input signal that maximizes the information retrieved from an experiment. Some of the initial contributions are discussed in \cite{Cox1958} and \cite{goopay76}. Since then, many contributions to the subject have been developed; see e.g. \cite{fedorov1972}, \cite{Whittle1973}, \cite{hildebrand2003}, \cite{gevers2005} and the references therein.

In this article, a new method for experiment design in nonlinear systems is presented, which extends the input design methods proposed in \cite{gopaluni2011} and \cite{valenzuela-rojas-hjalmarsson-13}. The objective is to design an experiment as a realization of a stationary process, such that the system is identified  with maximum accuracy as defined by a scalar function of the Fisher information matrix, and under the assumption that the input can adopt a finite set of values. %The method considers the optimization of a stationary probability mass function (pmf), from which the input signal is sampled.
 The assumption on the input class modifies the class of input sequences considered in \cite{gopaluni2011}. The optimization of the stationary probability mass function (pmf) is done by maximizing a scalar cost function of the information matrix over the feasible set of pmf's. Using concepts from graph theory \citep{zaman1983,johnson1975,tarjan1972}, we can express the feasible set of pmf's as a convex combination of the measures for the extreme points of the set. Therefore, the information matrix corresponding to a feasible pmf can be expressed as the convex combination of the information matrices associated with the extreme points of the feasible set. Since the exact computation of the information matrices for nonlinear systems is often intractable, we use particle methods to compute sampled information matrices for the extreme points of the feasible set. This allows us to extend the technique of \cite{valenzuela-rojas-hjalmarsson-13} to more general nonlinear model structures. An attractive property of the method is that the optimization problem is convex even for nonlinear systems. In addition, since the input is restricted to a finite set of possible values, the method can naturally handle amplitude limitations.

%As with most optimal input design methods, the one proposed in this contribution relies on knowledge of the true system. This difficulty can be overcome by implementing a robust experiment design scheme on top of it \citep{Rojas2007} or via an adaptive procedure, where the input signal is re-designed as more information is being collected from the system \citep{rojas2011adaptive}. Due to space limitations, however, we will not address these issues in the present article.

%In dynamical systems, the input design problem can be formulated as the optimization of a cost function related to the model application. The development of input design for dynamical systems has mostly been concerned with the linear case.
Previous results on input design have mostly been concerned with linear systems. A Markov chain approach to input design is presented in \cite{brighenti2009}, where the input is modelled as the output of a Markov chain. \cite{Suzuki2007} presents a time domain experiment design method for system identification. Linear matrix inequalities (LMI) are used to solve the input design problem in \cite{jansson2005} and \cite{lindqvist2000}. A robust approach for input design is presented in \cite{Rojas2007}, where the input signal is designed to optimize a cost function over a set where the true parameter is assumed to lie.%the feasible set of true parameters.

In recent years, the interest in input design for nonlinear systems has increased. The main problem here is that the frequency domain approach for experiment design used in linear systems is no longer valid. An analysis of input design for nonlinear systems using the knowledge of linear systems is considered in \cite{hjalmarsson2007}. In \cite{larsson2010} an input design method for a particular class of nonlinear systems is presented.  Input design for structured nonlinear systems is discussed in \cite{vincent2009}. \cite{gopaluni2011} introduces a particle filter method for input design in nonlinear systems. An analysis of input design for a class of Wiener systems is considered in \cite{decock-gevers-schoukens-13}. A graph theory approach for input design for output-error like nonlinear system is presented in \cite{valenzuela-rojas-hjalmarsson-13}. The results presented allow to design input signals when the system contains nonlinear functions, but the restrictions on the system dynamics and/or the input structure are the main limitations of most of the previous contributions. Moreover, with the exception of \cite{brighenti2009}, \cite{larsson2010} and \cite{valenzuela-rojas-hjalmarsson-13}, the proposed methods cannot handle amplitude limitations on the input signal, which could arise due to physical and/or safety reasons.

\section{Problem formulation}\label{sec:  2}
In this article, the objective is to design an input signal $u_{1:n_{\op{seq}}}:= \{ u_t\}_{t=1}^{n_{\op{seq}}}$, as a realization of a stationary process. This is done such that a state space model (SSM) can be identified with maximum accuracy as defined by a scalar function of the Fisher information matrix $\Ical _F$ \citep{ljung99}. An SSM with states $x_{1:T} := \{x_t\}_{t=1}^T$, inputs $u_{1:T}$ and measurements $y_{1:T}$ is given by
\begin{subequations}
\begin{align}
	x_{t}|x_{t-1}  &\sim  f_{\theta}(x_{t}|x_{t-1},u_{t-1}), \\
	y_{t}|x_t      &\sim  g_{\theta}(y_{t}|x_t,u_t).
\end{align}%
\label{eq:generalSSM}%
\end{subequations}%
Here, $f_{\theta}(\cdot)$ and $g_{\theta}(\cdot)$ denote known probability distributions parametrised by $\theta \in \Theta \subset \rb ^d$. For the remainder of this article, we make the rather restrictive albeit standard assumption that we know the initial state $x_0$ and the true model structure \eqref{eq:generalSSM} with true parameters $\theta _0$. Hence, we can write the joint distribution of states and measurements for \eqref{eq:generalSSM} as
\begin{align}
	p_{\theta}(x_{1:T},y_{1:T}|u_{1:T}) =
	\prod_{t=1}^{T} f_{\theta}(x_{t}|x_{t-1},u_{t-1})
    g_{\theta}(y_{t}|x_t,u_t).
	\label{eq:genSSMjointDist}
\end{align}
This quantity is used in the sequel for estimating $\Ical _F$ by
\begin{subequations}
\begin{align}
\label{eq: prob3}
\Ical _F &:= \Ebf _\theta \left\{ \mathcal{S}(\theta _0) \mathcal{S}^{\top}( \theta_0) \right\} \, , \\
\label{eq: prob2a}
\mathcal{S}(\theta _0) &:= \left. \nabla \, \textsf{log } l_\theta (y_{1:n_{\op{seq}}}) \right| _{\theta = \theta _0}\, ,
\end{align}%
\label{eq: fim}%
\end{subequations}%
\noindent where $l_\theta (y_{1:n_{\op{seq}}})$ and $\mathcal{S}(\theta)$ denote the likelihood function and the score function, respectively. Note, that the expected value in \eqref{eq: prob3} is with respect to the stochastic processes in \eqref{eq:generalSSM} and the realizations of $u_{1:n_{\op{seq}}}$. %and $\Ecal _{n_{\op{seq}}}:= \{e_{(n_{\op{seq}}-1)}, \, \ldots , \, e_0 \}$.
%In addition, the definition in \eqref{eq: fim} assumes that there exists a $\theta _0 \in \Theta$ such that the model \eqref{eq:generalSSM} becomes the true system when $\theta = \theta _0$ \citep{ljung99}, i.e., that there is no undermodelling. %; we will make this assumption in the sequel.

We note that \eqref{eq: prob3} depends on the cumulative density function (cdf) of $u_{1:n_{\op{seq}}}$, say $P_u(u_{1:n_{\op{seq}}})$. Therefore, %to obtain an optimal stationary input sequence $\Ucal _{n_{\op{seq}}} ^{\op{opt}}$,
 the input design problem is to find a cdf %cumulative distribution function
  $P^{\op{opt}}_u(u_{1:n_{\op{seq}}})$ which optimizes a scalar function of \eqref{eq: prob3}. We define this scalar function as $h_m : \, \rb ^{d \times d} \rightarrow \, \rb$. To obtain the desired results, $h_m$ must be a nondecreasing matrix function \cite[pp. 108]{boyvan04}. Different choices of $h_m$ have been proposed in the literature, see e.g. \cite{Rojas2007}; some examples are $h_m = \det$, and $h_m = -\op{tr} \{(\cdot)^{-1}\}$. In this work, we leave the selection of $h_m$ to the user. %Therefore, we can handle different optimality criteria by using the method introduced in this article.

Since $P^{\op{opt}}_u(u _{1:n_{\op{seq}}})$ has to be a stationary cdf, the optimization must be constrained to the set
\begin{multline}
\label{eq: prob5}
\Pcal := \left\{ P_u:\, \rb ^{n_{\op{seq}}} \rightarrow \rb | \, P_u(\xbf) \geq 0 , \, \forall \xbf \in \rb ^{n_{\op{seq}}}; \,  \right. \\
\left. P_u \text{ is monotone non-decreasing}\right. ; \, \\
\left. \lim _ {\substack{x_i \rightarrow \infty \\ i = \{1, \, \ldots , \, n_{\op{seq}}\}\\ \xbf = (x_1 , \, \ldots , \, x_{n_{\op{seq}}})}} P_u(\xbf) = 1 ;\,\right. \\
\left. \int _{v \in \rb} dP_u(v, \, \zbf) = \int _{v \in \rb} dP_u(\zbf , \, v) \, , \forall \zbf \in \rb ^{n_{\op{seq}}-1} \right\} \, .
\end{multline}
The last condition in \eqref{eq: prob5} (with slight abuse of notation) guarantees that $P_u \in \Pcal$ is the cdf of a stationary sequence \citep{zaman1983}.

To simplify our analysis, we will assume that $u_t$ can only adopt a finite number $c_{\op{seq}}$ of values. We define this set of values as $\Ccal$. With the previous assumption, we can define the following subset of $\Pcal$:% as
\begin{multline}
\label{eq: prob6}
\Pcal _\Ccal := \left\{ p_u:\, \Ccal ^{n_{\op{seq}}} \rightarrow \rb | \, p_u(\xbf) \geq 0 , \, \forall \xbf \in \Ccal ^{n_{\op{seq}}}; \,  \right. \\
\left. \sum _ {\xbf \in \Ccal ^{n_{\op{seq}}}} p_u(\xbf) = 1; \right.  \\
\left. \sum _{v \in \Ccal} p_u(v, \, \zbf) = \sum _{v \in \Ccal} p_u(\zbf , \, v) \, , \forall \zbf \in \Ccal ^{(n_{\op{seq}}-1)} \right\} \, .
\end{multline}
The set introduced in \eqref{eq: prob6} will constrain the pmf $p_u(u _{1:n_{\op{seq}}})$.

The problem described can be summarized as
\begin{prob}\label{prob: 1}
Design an optimal input signal $u_{1:n_{\op{seq}}} \in \Ccal ^{n_{\op{seq}}}$ as a realization from $p^{\op{opt}}_u(u_{1:n_{\op{seq}}})$, where% such that
\beq
\label{eq: prob7}
p^{\op{opt}}_u := \arg \max _{p_u \in \Pcal _\Ccal} h_m(\Ical _F(p_u)) \, ,
\eeq
with $h_m : \, \rb ^{d \times d} \rightarrow \, \rb$ a matrix nondecreasing function,
%\beq
%\label{eq: prob8}
%\Ical _F(p) = \dfrac{1}{\lambda _e} \sum _{u_{1:n_{\op{seq}}} \in \Ccal ^{n_{\op{seq}}}}  \sum_{t=1} ^{n_{\op{seq}}} \psi(\theta _0) \psi( \theta_0) ^T\, p(\Ucal _{n_{\op{seq}}}) \, ,
%\eeq
and $\Ical _F \in \rb ^{d \times d}$ defined as in \eqref{eq: fim}.
\fin
\end{prob}
%A solution for this problem will be discussed in the next section.
%%%%%%%%%%%%%%%%%%%%%%%%%%%%%%%%%%%%%%%%%%%%%%%%%%%%%%%%%%%%%%%%%%%%%%%%%%%%%%%%
\section{New input design method}\label{sec: 3}
In this section, we discuss the proposed input design method, which is based on three steps. In the first step, we calculate basis input signals, which are used to excite the system. In the second step, we iteratively calculate the information matrix estimate and the optimal weighting of the basis inputs in a Monte Carlo setting. In the third step, we generate an optimal input sequence using the estimated optimal weighting of the basis inputs. %introduce the extension to the input design methods presented in \citep{gopaluni2011,valenzuela-rojas-hjalmarsson-13}.
%%%%%%%%%%%%%%%%%%%%%%%%%%%%%%%%%%%%%%%%%%%%%%%%%%%%%%%%%%%%%%%%%%%%%%%%%%%%%%%%
%\subsection{Preliminaries on graph theory}\label{sec: 1a}
%%The results we will introduce in the coming sections require the knowledge of basic concepts in graph theory.
%The purpose in this subsection is to provide a brief background on graph theory to understand the discussion in the next subsections. The definitions presented here come from \cite[pp. 77]{johnson1975}.
%
%A \emph{directed graph} $\Gcal _{\Vcal} := (\Vcal ,\Xcal)$ consists of a nonempty and finite set of vertices (or nodes) $\Vcal$ and a set $\Xcal$ of ordered pairs of distinct vertices called \emph{edges}. A \emph{path} in $\Gcal _\Vcal$ is a sequence of vertices $p_{zu} := (z=z_1,\, z_2 ,\,  \ldots , \, z_k = u)$ such that $(z_i ,\, z_{i+1}) \in \Xcal$ for $i \in \{1, \, \ldots , \, k-1 \}$. A \emph{cycle} is a path in which the first and last vertices are identical. %A path is \emph{elementary} if no vertex appears twice.
%A cycle is elementary if no vertex but the first and last appears twice. Two elementary cycles are distinct if one is not a cyclic permutation of the other.
%%%%%%%%%%%%%%%%%%%%%%%%%%%%%%%%%%%%%%%%%%%%%%%%%%%%%%%%%%%%%%%%%%%%%%%%%%%%%%%%
\subsection{Graph theoretical input design}
Problem \ref{prob: 1} is often hard to solve explicitly since
\begin{itemize}
\item[(i)] we need to represent the elements in $\Pcal _\Ccal$ as a linear combination of its basis functions, and
\item[(ii)] the stationary pmf $p_u$ is of dimension $n_{\op{seq}}$, where $n_{\op{seq}}$ could potentially be very large.
\end{itemize}
These issues make Problem \ref{prob: 1} computationally intractable.

To solve issue (ii), %the dimensionality issue associated with $n_{\op{seq}}$,
 we assume that $p_u$ is an extension from the subspace of stationary pmf's of memory length $n_m$, where $n_m << n_{\op{seq}}$. %This assumption makes Problem \ref{prob: 1} feasible to solve.

To address issue (i), notice that all the elements in $\Pcal _\Ccal$  can be represented as a convex combination of its extreme points \citep{valenzuela-rojas-hjalmarsson-13}. %is a convex set. In particular, $\Pcal _\Ccal$ is a polyhedron \cite[pp. 31]{boyvan04}. %Since $\Pcal _\Ccal$ is a polyhedron, then
 %Hence, any element of $\Pcal _\Ccal$ can be described as a convex combination of the extreme points of $\Pcal _\Ccal$ \cite[pp. 24]{boyvan04}. Therefore, if we define $\Vcal_{\Pcal _\Ccal}$ as the set of all the extreme points of $\Pcal _\Ccal$, composed by $n_\vcal$ elements, then for all $f_u \in \Pcal _\Ccal$ we have
%\beq
%\label{eq: prob9}
%f_u = \sum _{i=1}^{n_{\vcal}} \alpha _i \, v_i \, ,
%\eeq
%where $\alpha _i \geq 0$, $i \in \{1, \ldots , \, n_\vcal\}$, %and
%\beq
%\label{eq: prob10}
%\sum _{i=1}^{n_{\vcal}} \alpha _i = 1 \, ,
%\eeq
%and $v_i  \in \Vcal_{\Pcal _\Ccal}$, for all $i \in \{1, \ldots , \, n_\vcal\}$.
 %In the following,
  We will refer to $\Vcal_{\Pcal _\Ccal} := \{ v_1 , \, \ldots , \, v_{n_\vcal}\}$ as the set of the extreme points of $\Pcal _\Ccal$. %, which is composed by

%Equation \eqref{eq: prob9} says that all the elements in $\Pcal _\Ccal$ can be described by using $n_{\vcal}$ elements in the set $\Vcal_{\Pcal _\Ccal}$. %Therefore, if we are able to find all the elements in the set $\Vcal_{\Pcal _\Ccal}$, then we can overcome drawback 1), which reduces significantly the computational cost of solving Problem 1.

To find all the elements in $\Vcal_{\Pcal _\Ccal}$, we will make use of graph theory %. Indeed, we can analyze the set $\Ccal ^{n_{m}}$
 as follows. $\Ccal ^{n_{m}}$ is composed of $(c_{\op{seq}} )^{n_{m}}$ elements. Each element in $\Ccal ^{n_{m}}$ can be viewed as one node in a graph. In addition, the transitions (edges) between the elements in $\Ccal ^{n_{m}}$ are given by the feasible values of $u_{k+1}$ when we move from $(u_{k-n_m+1} , \ldots , \, u_k)$ to $(u_{k-n_m+2} , \ldots , \, u_{k+1})$, for all integers $k \geq 0$. %Since the transitions among the elements in $\Ccal ^{n_{\op{seq}}}$ are clearly defined,
 %The edges between the elements in $\Ccal ^{n_{m}}$ denote the possible transitions between the states, represented by the nodes of the graph. %must be consistent with this assumption.
  Figure \ref{fig: 2} illustrates this idea, when $c_{\op{seq}}= 2$, $n_{m}=2$, and $\Ccal = \{0,\, 1\}$. From this figure we can see that, if we are in node $(0,\, 1)$ at time $t$, then we can only transit to node $(1, \, 0)$ or $(1, \, 1)$ at time $t+1$. %If the last one is not satisfied, then the graph is not derived from $\Ccal ^{n_{\op{seq}}}$. %a set with sequences.
   \begin{figure}[t]%hpb]
      \centering
      \input{graph.pstex_t}
      \caption{Example of graph derived from $\Ccal ^{n_{m}}$, with $c_{\op{seq}}=2$, $n_{m}=2$, and $\Ccal := \{0, \, 1\}$.} %(we omit self-edges in $(0,\, 0)$ and $(1, \, 1)$ to improve legibility).}
      \label{fig: 2}
   \end{figure}

%The idea to use graph theory
To find all the elements in $\Vcal_{\Pcal _\Ccal}$ we rely on the concept of prime cycles. A \emph{prime cycle} is an elementary cycle whose set of nodes do not have a proper subset which is an elementary cycle \cite[pp. 678]{zaman1983}. It has been proved that the prime cycles of a graph describe all the elements in the set $\Vcal_{\Pcal _\Ccal}$ \cite[Theorem 6]{zaman1983}. In other words, each prime cycle defines one element $v_j \in \Vcal_{\Pcal _\Ccal}$. Furthermore, each $v_j$ corresponds to a uniform distribution whose support is the set of elements of its prime cycle, for all $j \in \{1, \ldots , \, n_\vcal\}$ \cite[pp. 681]{zaman1983}. Therefore, the elements in $\Vcal_{\Pcal _\Ccal}$ can be described by finding all the prime cycles associated with the stationary graph $\Gcal _{\Ccal ^{n_{m}}}$ drawn from $\Ccal ^{n_{m}}$.

\begin{figure}[t]%hpb]
      \centering
      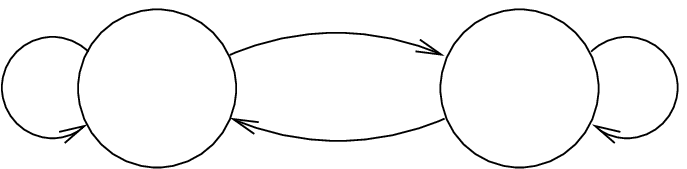
      \caption{Example of graph derived from $\Ccal ^{n_{m}}$, with $c_{\op{seq}}=2$, $n_{m}=1$, and $\Ccal := \{0, \, 1\}$.}
      \label{fig: 2a}
   \end{figure}

It is known that all the prime cycles associated with $\Gcal _{\Ccal ^{n_{m}}}$ can be derived from the elementary cycles associated with $\Gcal _{\Ccal ^{(n_{m}-1)}}$ \cite[Lemma 4]{zaman1983}, which can be found by using existing algorithms\footnote{For the examples in Section \ref{sec: 4}, we have used the algorithm presented in \cite[pp. 79--80]{johnson1975} complemented with the one proposed by \cite[pp. 157]{tarjan1972}.}. %An \emph{elementary circuit} is a path where no vertex but the first and last appears twice \cite[pp. 77]{johnson1975}.
 %In the literature there are many algorithms for finding all the elementary cycles in a graph\footnote{For the examples in Section \ref{sec: 4}, we have used the algorithm presented in \cite[pp. 79--80]{johnson1975} complemented with the one proposed by \cite[pp. 157]{tarjan1972}.}.
%Once all the elementary cycles of $\Gcal _{\Ccal ^{(n_{m}-1)}}$ are found, we can find all the prime cycles associated to $\Gcal _{\Ccal ^{n_{m}}}$ according to \cite[Lemma 4]{zaman1983}.
 To illustrate this, we consider the graph depicted in Figure \ref{fig: 2a}. One elementary cycle for this graph is given by $(0,\, 1, \, 0)$. %Therefore
 Using \cite[Lemma 4]{zaman1983}, the elements of one prime cycle for the graph $\Gcal _{\Ccal ^{2}}$ are obtained
 as a concatenation of the elements in the elementary cycle $(0,\, 1, \, 0)$. Hence, the prime cycle in $\Gcal _{\Ccal ^{2}}$ associated with this elementary cycle is given by $((0,\, 1), \, (1,\, 0), \, (0,\, 1))$.

%With all the prime cycles clearly defined for $\Gcal _{\Ccal ^{n_{m}}}$, then all the elements in the set $\Vcal_{\Pcal_\Ccal}$ are found. Hence, %Since we know all the elements in $\Vcal_{\Pcal_\Ccal}$,
%we can use \eqref{eq: prob9} to describe all the elements in $\Pcal _\Ccal$. Thus, the solution described here presents a computationally feasible method to address the first issue.
%
%%With the method presented to solve drawback 1), we can also solve drawback 2).
%Since we know the distribution $v_i$ for each prime cycle, with $v_i  \in \Vcal _{\Pcal _\Ccal}$, we can generate an input signal $\{u_t ^i\}_{t=0} ^{t=N}$ drawn from $v_i $, so that
%\begin{align}
%\nonumber
%\Ical _F ^{(i)} &:= \dfrac{1}{\lambda _e} \sum _{\Ucal _{n_{m}} \in \Ccal ^{n_{m}}} \sum_{t=1}^{n_m} \psi(\theta _0) \psi( \theta_0) ^T\, v_i(\Ucal _{n_{m}}) \\
%\label{eq: prob11}
%&\approx \dfrac{1}{\lambda _e  \, N} \sum _{t=1} ^{N} \psi(\theta _0) \psi( \theta_0) ^T \, ,
%\end{align}
%for all $i \in \{1, \ldots , \, n_\Vcal\}$, and $N$ sufficiently large\footnote{Note that $N$ is the number of Monte Carlo simulations to compute \eqref{eq: prob11}, and it is not necessarily equal to the length of the experiment $n_{\op{seq}}$.} (in relation to the length of the prime cycles). Notice that $\psi(\theta _0)$ depends implicitly on $\{u_t ^i\}_{t=0} ^{t=N}$ through \eqref{eq: prob2a}-\eqref{eq: prob2b}. Furthermore, each $\Ical _F ^{(i)}$ is associated to the $i$-th prime cycle, for all $i \in \{1, \ldots , \, n_\Vcal\}$.

Since we know the prime cycles, it is possible to generate an input sequence $\{u_t ^j\}_{t=0} ^{T}$ from $v_j$, which will be referred to as the \emph{basis inputs}. As an example, we use the graph depicted in Figure \ref{fig: 2}. One prime cycle for this graph is given by $((0,\, 1), \, (1,\, 0), \, (0, \, 1))$. Therefore, the sequence $\{u_t ^j\}_{t=0} ^{T}$ is given by taking the last element of each node, i.e., $\{u_t ^j\}_{t=0} ^{T} = \{1,\, 0 , \, 1, \, 0 , \, \ldots, \,  ((-1)^{T}+1)/2\}$.

Given $\{ u_t ^j\}_{t=0}^{T}$, %the input sequences associated with each $v_j \in \Vcal_{\Pcal_\Ccal}$,
 we can use them to obtain the corresponding information matrix for $v_j \in \Vcal_{\Pcal_\Ccal}$, say $\Ical _F ^{(j)}$. However, in general the matrix $\Ical _F ^{(j)}$ cannot be computed explicitly. %Indeed, the expression \eqref{eq: prob2a} is often intractable.
 To overcome this problem, we use Sequential Monte Carlo methods to approximate $\Ical _F ^{(j)}$, as discussed in the next subsection.
%The approximation of each $\Ical _F ^{(i)}$ given by \eqref{eq: prob11} reduces the sum \eqref{eq: prob8} from dimension $n_{\op{seq}}$ to dimension 1. This simplification reduces significantly the computation effort to obtain \eqref{eq: prob8}. With this approach, issue 2) is also addressed.
%%%%%%%%%%%%%%%%%%%%%%%%%%%%%%%%%%%%%%%%%%%%%%%%%%%%%%%%%%%%%%%%%%%%%%%%%%%%%%%%
% SMC part
%\input{smcsection.tex}
\subsection{Estimation of the Score function}
%\item Introduce SMC methods as a tool to estimate the latent states of a nonlinear state space model. This is done by simulation of the state dynamics and by comparing the simulated states (called particle) with the given measurements.
Sequential Monte Carlo (SMC) methods are a family of methods that can be used e.g.\ to estimate the filtering and smoothing distributions in SSMs. General introductions to SMC samplers are given in e.g.\
 \cite{DoucetJohansen2011} and \cite{DelMoralDoucetJasra2006}. Here, we introduce the \textit{auxiliary particle filter} (APF) \citep{PittShephard1999} and the fixed-lag (FL) particle smoother \citep{KitagawaSato2001} to estimate the score function for \eqref{eq:generalSSM}. In the next subsection, the score function estimates are used with \eqref{eq: fim} to estimate the information matrix.

The APF estimates the smoothing distribution by
\begin{align}
	\widehat{p}_{\theta}(\dn x_{1:t} | y_{1:t})
	:=
	\sum_{i=1}^N
	\frac{w_{t|t}^{(i)}}{\sum_{k=1}^N w_{t|t}^{(k)}}
	\delta_{x_{1:t}^{(i)}}(\dn x_{1:t}),
	\label{eq:APFfilteringDist}
\end{align}
where the weights $w_{t|t}^{(i)}$ and the particle trajectories $x^{(i)}_{1:t}$ are computed by the APF as a article system, $\{x_{1:t}^{(i)},w_{t|t}^{(i)}\}_{i=1}^N$. Here, $\delta_{z}(\dn x_{1:t})$ denotes the Dirac measure at $z$.

The particle system is sequentially computed using two steps: (i) sampling/propagation and (ii) weighting. The first step can be seen as sampling from a proposal kernel,
\begin{align}
	\{a_t^{(i)}, x^{(i)}_t\}
	&\sim
	\frac
	{w^{a_t}_{t-1|t-1}}
	{\sum_{k=1}^N w^{(k)}_{t-1|t-1}}
	R_{\theta,t}(x_t | x^{a_t}_{t-1}, u_{t-1}), \label{eq:APFmutation}
\end{align}
where we append the sampled particle to the trajectory by $x_{1:t}^{(i)}=\{x_{1:t-1}^{(i)},x_t^{(i)}\}$. Here, $R_{\theta,t}(\cdot)$ denotes the propagation kernel and the \textit{ancestor index} $a_t^{(i)}$ denotes the index of the \textit{ancestor} at time $t-1$ of particle $x_t ^{(i)}$. In the second step, we calculate the (unnormalised) importance weights, % are given by a weighting function
\begin{align}
	w_{t|t}^{(i)}
	\triangleq
	\frac
	{ g_{\theta}(y_t|x_t^{(i)},u_t) f_{\theta}(x_t|x_{t-1}^{(i)},u_{t-1}) }
	{ R_{\theta,t}(x_t|x_{t-1}^{(i)},u_{t-1}) }\, .
	\label{eq:APFweights}
\end{align}
%Repeating, these two steps generates the particle system sequentially over time.

%\item The score can be estimated using Fisher's identity which requires an estimate of the marginal smoothing distribution.
SMC methods can be used to compute an estimate of the score function in combination with \textit{Fisher's identity} \citep{Fisher1925,CappeMoulinesRyden2005},
\begin{align*}
	\nabla \textsf{log } l_{\theta}(y_{1:T}) =
	\mathbf{E}_{\theta}[	\nabla \textsf{log } p_{\theta}(x_{1:T},y_{1:T}|u_{1:T}) | y_{1:T},u_{1:T}].
\end{align*}
Inserting \eqref{eq:genSSMjointDist}, we obtain
\begin{align*}
\nabla \textsf{log } l_{\theta}(y_{1:T}) =
\sum_{t=1}^{T} \dint
\nabla \textsf{log } g_{\theta}(y_{t}|x_t,u_t)
p_{\theta}(x_t|y_{1:T}) \dd x_t\\
+
\sum_{t=1}^{T} \dint
\nabla \textsf{log } f_{\theta}(x_{t}|x_{t-1},u_{t-1})
p_{\theta}(x_{t-1:t}|y_{1:T}) \dd x_{t-1:t},
\end{align*}
which depends on the one-step and two-step marginal smoothing densities. The APF can be used the estimate these quantities but this leads to poor estimates with high variance, due to problems with \textit{particle degeneracy}.

%\item For this we use the FL-smoother (see Dahlin, Lindsten and Schön (2013b) for the details) to estimate the score functjointion. Give the expressions.
Instead, we use an FL-smoother to estimate the smoothing densities, which reduces the variance of the score estimates \citep{OlssonCappeDoucMoulines2008}. The fixed-lag smoother assumes that
\begin{align*}
p_{\theta}(x_t|y_{1:T},u_{1:T}) \approx p_{\theta}(x_t|y_{1:\kappa_t},u_{1:\kappa_t}),
\end{align*}
for $\kappa_t = \min(t+\Delta,T)$ with some fixed-lag $\Delta$. This means that measurements after some time has a negligible effect on the state, see \citep{DahlinLindstenSchon2013b} for more details about the FL-smoother and its use for score estimation. The resulting expression is obtained as
\begin{align}
	\widehat{\mathcal{S}}(\theta)
	&:=
	\sum_{t=1}^{T}
	\left[
	\mathcal{L}_{1,t}(\theta)
	+
	\mathcal{L}_{2,t}(\theta)
	\right], \text{ with}
	\label{eq:FisherScoreParticleApproximation}
	\\
	\mathcal{L}_{1,t}(\theta)
	&:= \sum_{i=1}^N w_{\kappa_t|\kappa_t}^{(i)}
	\nabla \textsf{log } g_{\theta}(y_{t}|x_{t}^{a_{k,t}^{(i)}},u_t),
	\nonumber
	\\	
	\mathcal{L}_{2,t}(\theta)
	&:= \sum_{i=1}^N
	w_{\kappa_t|\kappa_t}^{(i)} \nabla \textsf{log } f_{\theta}(x_{t}^{a_{k,t}^{(i)}}|x_{t-1}^{a_{k,t-1}^{(i)}},u_{t-1}),
	\nonumber
\end{align}%
where, $a_{k,t}^{(i)}$ denotes the particle at time $t$ which is the ancestor of particle $i$ at time $k$. The complete procedure for estimating the score function using the FL smoother is outlined in Algorithm \ref{alg:ScoreEst}.

%\item Write out the explicit algorithm for score estimation.
\begin{algorithm}[!t]
\caption{\textsf{Score estimation using particle FL-smoothing}}
\begin{small}
\textsf{Input:} The SSM on the form \eqref{eq:generalSSM} with measurements $y_{1:T}$ and inputs $u_{1:T}$. The propagation kernel $R_{\theta,t}(\cdot)$ and the number of particles $N$. \\
\textsf{Output:} An estimate of the score function $\widehat{\mathcal{S}}(\theta)$.
\algrule[.4pt]
\begin{itemize}[leftmargin=*]
	\item \textsf{Run the auxiliary particle filter}
	\item[] Initialise particles $x_0^{(i)}$ for $i=1,\ldots,N$.
	\item[] \textbf{for} $t=1,\ldots,T$ \textbf{do}
	\begin{itemize}
		\item[-] Sample and propagate the particles using \eqref{eq:APFmutation}.
		\item[-] Calculate the particle weights using \eqref{eq:APFweights}.
	\end{itemize}
	\item[] \textbf{end for}
	\item \textsf{Run the fixed-lag particle smoother}
	\item[] \textbf{for} $t=1,\ldots,T$ \textbf{do}
	\begin{itemize}
		\item[-] \textbf{set} $\kappa_t = \min\{T,t+\Delta\}$.
		\item[-] Recover the ancestor indices $a_{\kappa_t,t}$.
	\end{itemize}	
	\item[] \textbf{end for}
	\item Compute the score function estimate using \eqref{eq:FisherScoreParticleApproximation}.
\end{itemize}
\end{small}
\label{alg:ScoreEst}
\end{algorithm}

\subsection{Monte Carlo-based optimisation}
Given $\{\Ical _F ^{(j)}\}_{j=1}^{n_\Vcal}$ associated with the elements in $\Vcal_{\Pcal _\Ccal}$, we can find the corresponding information matrix associated with any element in $\Pcal _\Ccal$ as a convex combination of the $\Ical _F ^{(j)}$'s. By defining $\gamma := \{\alpha _1 , \ldots , \, \alpha _{n_\Vcal}\} \in \rb ^{n_{\Vcal}}$, we introduce $\Ical _F ^{\op{app}}(\gamma)$ as the information matrix associated with one element in $\Pcal _\Ccal$ for a given $\gamma$ such that $\alpha_j \geq 0$, $j \in \{1, \, \ldots , \, n_{\Vcal}\}$, {$\sum _{j=1} ^{n_{\Vcal}} \alpha _j = 1$}. Therefore, finding the optimal $\Ical _F ^{\op{app}}(\gamma)$ is equivalent to determining the optimal weighting vector $\gamma$.

Hence, we can rewrite Problem~\ref{prob: 1} as %the following problem:
\begin{subequations}
\begin{align}
\label{eq: prob12aa}
\gamma ^{\op{opt}} &= \arg \max _{\gamma \in \rb ^{n_{\Vcal}}} h_m(\Ical _F ^{\op{app}}(\gamma)) \, , \\
\label{eq: prob12a}
\text{st.} \quad &\Ical _F ^{\op{app}}(\gamma) :=  \sum _{j=1} ^{n_{\Vcal}} \alpha _j \, \Ical _F ^{(j)} \, , \\
&\sum _{j=1} ^{n_{\Vcal}} \alpha _j = 1 \, , \\
\label{eq: prob13}
&\alpha _j \geq 0 \, , \text{ for all } j \in \{1, \ldots , \, n_{\Vcal}\} \, ,
\end{align}%
\label{eq: prob12}%
\end{subequations}%
\noindent To solve the optimisation problem \eqref{eq: prob12}, we need to estimate the information matrix for each basis input.

In the SMC literature, the observed information matrix is often estimated by the use of \textit{Louis' identity} \citep{Louis1982,CappeMoulinesRyden2005}. However, this approach does not guarantee that the information matrix estimate is positive semi-definite. In the authors' experience, this standard approach also leads to poor accuracy in the estimates.

Instead, we make use of the fact that the information matrix can be expressed as \eqref{eq: fim}, i.e.\ the variance of the score function. Hence, a straight-forward method for estimating the information matrix is to use the Monte Carlo covariance estimator over some realisations of the system. If we denote each Monte Carlo estimate of the score function by $\widehat{\mathcal{S}}_m(\theta)$, the information matrix can be estimated using
\begin{align}
\label{eq: var_score_app}
	\widehat{\mathcal{I}}_F = \frac{1}{M-1} \sum_{m=1}^M \widehat{\mathcal{S}}_m(\theta) \widehat{\mathcal{S}}^{\top}_m(\theta),
\end{align}
where $M$ denotes the number of score estimates. Note, that this is an estimate of Fisher's information matrix as the Monte Carlo estimator averages over the system realisations. The estimate is positive semi-definite by construction but inherits some bias from the FL-smoother, see \cite{OlssonCappeDoucMoulines2008} for more information. This problem can be handled by using more computationally costly particle smoother. Later, we present results indicating that this bias does not effect the resulting input signal to any large extent.

The information matrix estimate in \eqref{eq: var_score_app} can be used to estimate $\Ical _F ^{(j)}$ for each basis input. A simple solution is therefore to plug-in the estimates and solve the convex optimisation problem \eqref{eq: prob12} using some standard solver. However, by doing this we neglect the stochastic nature of the estimates and disregard the uncertainty. In practice, this leads to bad estimates of $\gamma$.

Instead, we propose the use of a Monte Carlo method which iterates two different steps over $K$ iterations. In step (a), we compute the information matrix estimates $\Ical _F ^{(j)}$ for each input using \eqref{eq: var_score_app}. In step (b), we solve the optimisation problem in \eqref{eq: prob12} using the estimates to obtain $\gamma_k$ at iteration $k$. The estimate of the optimal weighting vector $\gamma^{\star}$ is found using the sample mean of $\gamma=\{\gamma_1,\ldots,\gamma_K\}$, which can be complemented with an confidence interval (CI). Such CI could be useful in determining which of the basis inputs that are significant and should be included in the optimal input sequence. The outline of the complete procedure is presented in Algorithm \ref{alg:MCopt}.

\begin{algorithm}[!t]
\caption{\textsf{Optimal input estimating using Monte Carlo}}
\begin{small}
\textsf{Input:} The inputs for Algorithm \ref{alg:ScoreEst}, the number of Monte Carlo runs $K$ and the size of each batch $M$.\\
\textsf{Output:} An estimate of the optimal weighting $\gamma^{\star}$ of the basis inputs.
\algrule[.4pt]
\begin{itemize}[leftmargin=*]
	\item \textbf{for} $k=1,\ldots,K$ \textbf{do}
	\begin{itemize}
	\item[-] Generate $M$ samples of the score function using Algorithm \ref{alg:ScoreEst} for each basis input.
	\item[-] Estimate the information matrix by \eqref{eq: var_score_app} for each basis input.
	\item[-] Solve the problem in \eqref{eq: prob12}.
	\item[-] \textbf{set} $\gamma_k$ as the weighting factors obtained from the solver.
	\end{itemize}
	\item[] \textbf{end for}
	\item Compute the sample mean of $\gamma=\{\gamma_1,\ldots,\gamma_K\}$, denote it as $\gamma^{\star}$.
\end{itemize}
\end{small}
\label{alg:MCopt}
\end{algorithm}
%%%%%%%%%%%%%%%%%%%%%%%%%%%%%%%%%%%%%%%%%%%%%%%%%%%%%%%%%%%%%%%%%%%%%%%%%%%%%%%
\subsection{Summary of the method}
The proposed method for designing of input signals in $\Ccal ^{n_{m}}$ is summarized in Algorithm~\ref{alg:inputmethod}.
%\begin{algorithm}[!t]
%\caption{\textsf{New input design method}}
%\begin{small}
%\textsf{Input:} The values for the input $\Ccal$, the memory $n_m$, the number of input samples $N_{\op{sim}}$, the number of particles $N$, the number of score samples $M$, and the number of Monte Carlo runs $K$. \\
%\textsf{Output:} The optimal weights $\gamma ^{\star}$. %An estimate of the likelihood function $\LfuncEst[]$.
%\algrule[.4pt]
%    \begin{enumerate}
%    \item Compute all the elementary cycles of $\Gcal _{\Ccal ^{(n_{m}-1)}}$ by using, e.g., \cite[pp. 79--80]{johnson1975}, \cite[pp. 157]{tarjan1972}.
%    \item Compute all the prime cycles of $\Gcal _{\Ccal ^{n_{m}}}$ from the elementary cycles of $\Gcal _{\Ccal ^{(n_{m}-1)}}$ as explained above (c.f. \cite[Lemma 4]{zaman1983}).
%    \item Generate the input signals $\{u_t ^j\}_{t=0} ^{N_{\op{sim}}}$ from the prime cycles of $\Gcal _{\Ccal ^{n_{m}}}$, for each $j \in \{1, \, \ldots, \, n_\Vcal\}$.
%    %\item For each $i \in \{1, \, \ldots, \, n_\Vcal\}$, approximate $\Ical _F ^{(i)}$ using \eqref{eq: var_score_app}.
%    \item Execute Algorithm~\ref{alg:MCopt}.
%    \end{enumerate}
%\end{small}
%\label{alg:inputmethod}
%\end{algorithm}
\begin{algorithm}[!t]
\caption{\textsf{New input design method}}
\begin{small}
\textsf{Input:} The values for the input $\Ccal$, the memory $n_m$ and the number of input samples $T$. The inputs to Algorithm \ref{alg:MCopt}. \\
\textsf{Output:} An estimate of the optimal weighting $\gamma^{\star}$ of the basis inputs.
\algrule[.4pt]
    \begin{itemize}[leftmargin=*]
    \item Compute all the elementary cycles of $\Gcal _{\Ccal ^{(n_{m}-1)}}$ by using, e.g., \cite[pp. 79--80]{johnson1975}, \cite[pp. 157]{tarjan1972}.
    \item Compute all the prime cycles of $\Gcal _{\Ccal ^{n_{m}}}$ from the elementary cycles of $\Gcal _{\Ccal ^{(n_{m}-1)}}$ as explained above (c.f. \cite[Lemma 4]{zaman1983}).
    \item Generate the input signals $\{u_t ^j\}_{t=0} ^{T}$ from the prime cycles of $\Gcal _{\Ccal ^{n_{m}}}$, for each $j \in \{1, \, \ldots, \, n_\Vcal\}$.
    %\item For each $i \in \{1, \, \ldots, \, n_\Vcal\}$, approximate $\Ical _F ^{(i)}$ using \eqref{eq: var_score_app}.
    \item Execute Algorithm~\ref{alg:MCopt}.
    \end{itemize}
\end{small}
\label{alg:inputmethod}
\end{algorithm}
The algorithm computes $\gamma ^{\star}$ which defines the optimal pmf $p^{\op{opt}}_u(u _{1: n_{m}})$ as a convex combination of the measures associated with the elements in $\Vcal _{\Pcal _\Ccal}$, with $\gamma ^{\star}$ as the weighting vector. %via \eqref{eq: prob9}.
 Notice that $\Ical _F ^{\op{app}}(\gamma)$ in \eqref{eq: prob12a} is linear in the decision variables. Therefore, the optimization \eqref{eq: prob12} is convex.

\section{Numerical examples}\label{sec: 4}
The following examples present some applications of the proposed input design method.
\begin{exmp} \label{ex: 1}
Consider the linear Gaussian state space (LGSS) system with parameters $\theta=\{\theta_1,\theta_2\}$,
\begin{align*}
x_{t+1} &= \theta_1 x_t + u_t + v_t, &v_t \sim \mathcal{N}(0,\theta_2^2), \\
y_t &= x_t +e_t, &e_t \sim \mathcal{N}(0,0.1^2),
\end{align*}
where the true parameters are $\theta_0=\{0.5, 0.1\}$. We design experiments to identify $\theta$ with $n_{\op{seq}} = 5 \cdot 10^3$ time steps, memory length $n_m = 2$, and an input assuming values $\Ccal = \{-1, \, 0, \, 1\}$. The optimal experiments maximize $h_m(\Ical_F ^{\op{app}}(\gamma)) = \textsf{det} (\Ical_F ^{\op{app}} (\gamma))$, and $h_m(\Ical_F ^{\op{app}}(\gamma)) = -\textsf{tr} \left\{ (\Ical_F ^{\op{app}} (\gamma) )^{-1} \right\}$.

We generate $\{u_t ^j\}_{t=0} ^{T}$ for each $v_j \in \Vcal_{\Pcal_\Ccal}$ ($T = 10^2$) to compute the approximation \eqref{eq: var_score_app} for each $\Ical _F ^{(j)}$. Finally, the optimal input $u_{1:n_{\op{seq}}}$ is computed by running a Markov chain with $p^{\op{opt}} _u(u _{1:n_{m}})$ as stationary pmf, %\eqref{eq: prob14},
 where we discard the first $2 \cdot 10^6$ samples and keep the last $n_{\op{seq}} = 5 \cdot 10^3$ ones. In addition, we consider $K=100$, $M = 5\cdot 10^3$ and $N = 10^3$.
\begin{table}[t]
\caption{$h_m(\widehat{\Ical}_F )$, % for different input sequences,
 Example \ref{ex: 1}.}
\centering
\begin{tabular}{c|cc}
  Input / $h_m(\widehat{\Ical}_F )$ & $\textsf{log } \textsf{det} (\widehat{\Ical}_F )$ & $\textsf{tr} \left\{ (\widehat{\Ical}_F  )^{-1} \right\}$ \\
  \hline
  Optimal (det) & $20.67(0.01)$ & $1.51 \cdot 10^{-4} (5.18 \cdot 10^{-7})$ \\
  Optimal (tr)  & $20.82(0.01)$ & $1.32 \cdot 10^{-4} (4.45 \cdot 10^{-7})$ \\
  Binary        & $\mathbf{20.91}(0.01)$ & $\mathbf{1.21 \cdot 10^{-4}} (4.51 \cdot 10^{-7})$ \\
  Uniform       & $19.38(0.01)$ & $5.32 \cdot 10^{-4} (2.12 \cdot 10^{-6})$ \\
\end{tabular}
\label{tab: inputs}
\end{table}
As a benchmark, we generate $n_{\op{seq}}$ input samples from uniformly distributed white noise with support $[-1, \, 1]$, and the same amount of samples from binary white noise with values $\{-1, \, 1\}$. These input samples are employed to compute an approximation of $\Ical _F$ via \eqref{eq: var_score_app}.

Table \ref{tab: inputs} condenses the results obtained for each input sequence, where \emph{Optimal (det)} and \emph{Optimal (tr)} represent the results for the input sequences obtained from optimizing $\textsf{det} (\Ical_F ^{\op{app}} (\gamma))$, and  $-\textsf{tr} \left\{ (\Ical_F ^{\op{app}} (\gamma) )^{-1} \right\}$, respectively. The $95\%$ confidence intervals are given as $\pm$ the value in the parentheses. From the data we conclude that, for this particular example, the binary white noise seems to be the best input sequence. Indeed, the proposed input design method tries to mimic the binary white noise excitation, which is clear from the numbers in Table \ref{tab: inputs}.
\fin
\end{exmp}
%The previous example showed that the extended input design method recovers a binary white noise sequence when it is the optimal experiment. The next example shows an application of the input design method in nonlinear systems.
\begin{exmp}\label{ex: 2}
In this example we consider the system in \cite[Section 6]{gopaluni2011}, given by
\begin{align*}
x_{t+1} &= \theta_1 x_t +\dfrac{x_t}{ \theta_2 + x_t^2} + u_t + v_t, & v_t \sim \mathcal{N}(0,0.1^2),  \\
y_t &= \frac{1}{2} x_t + \frac{2}{5} x_t^2 + e_t, & e_t \sim \mathcal{N}(0,0.1^2),
\end{align*}
where $\theta = \{\theta_1,\theta_2\}$ denotes the parameters with true values $\theta_0 = \{ 0.7,0.6\}$. We design an experiment with the same settings as in Example~\ref{ex: 1} maximizing $h_m(\Ical_F ^{\op{app}}(\gamma)) = \textsf{det} (\Ical_F ^{\op{app}} (\gamma))$. A typical input realization obtained from the proposed input design method is presented in Figure~\ref{fig: input_det}.
%In this example we consider the system introduced in \cite[Section 6]{gopaluni2011}, given by
%\begin{subequations}
%\begin{align}
%\label{eq: growth_state}
%x_{t+1} &= a\, x_t +\dfrac{x_t}{b +x_t ^2} +u_t +e_t \, ,  \\
%\label{eq: growth_output}
%y_t &= c\, x_t +d\, x_t^2 +n_t \, ,
%\end{align}
%\label{eq:growth}
%\end{subequations}
%where $\theta _0^T = \begin{bmatrix} a & b & c & d\end{bmatrix} = \begin{bmatrix} 0.7 & 0.6 & 0.5 & 0.4\end{bmatrix}$ are the unknown parameters. Here, $\{e_t\}$ and $\{n_t\}$ are white Gaussian stochastic processes with zero mean and variances $\lambda _e = \lambda _n = 0.01$. We design experiments to identify $\theta_0$, with $n_{\op{seq}} = 5\cdot 10^3$, $n_m = 2$, and $\Ccal = \{-1, \, 0, \, 1\}$. For this example, the optimal experiment maximizes $h_m(\Ical_F ^{\op{app}}(\gamma)) = \textsf{det} (\Ical_F ^{\op{app}} (\gamma))$. %, and $h_m(\Ical_F ^{\op{app}}(\gamma)) = -\textsf{tr} \left\{ (\Ical_F ^{\op{app}} (\gamma) )^{-1} \right\}$.
\begin{figure}[t]
{\centering
\includegraphics[width = 0.45\textwidth]{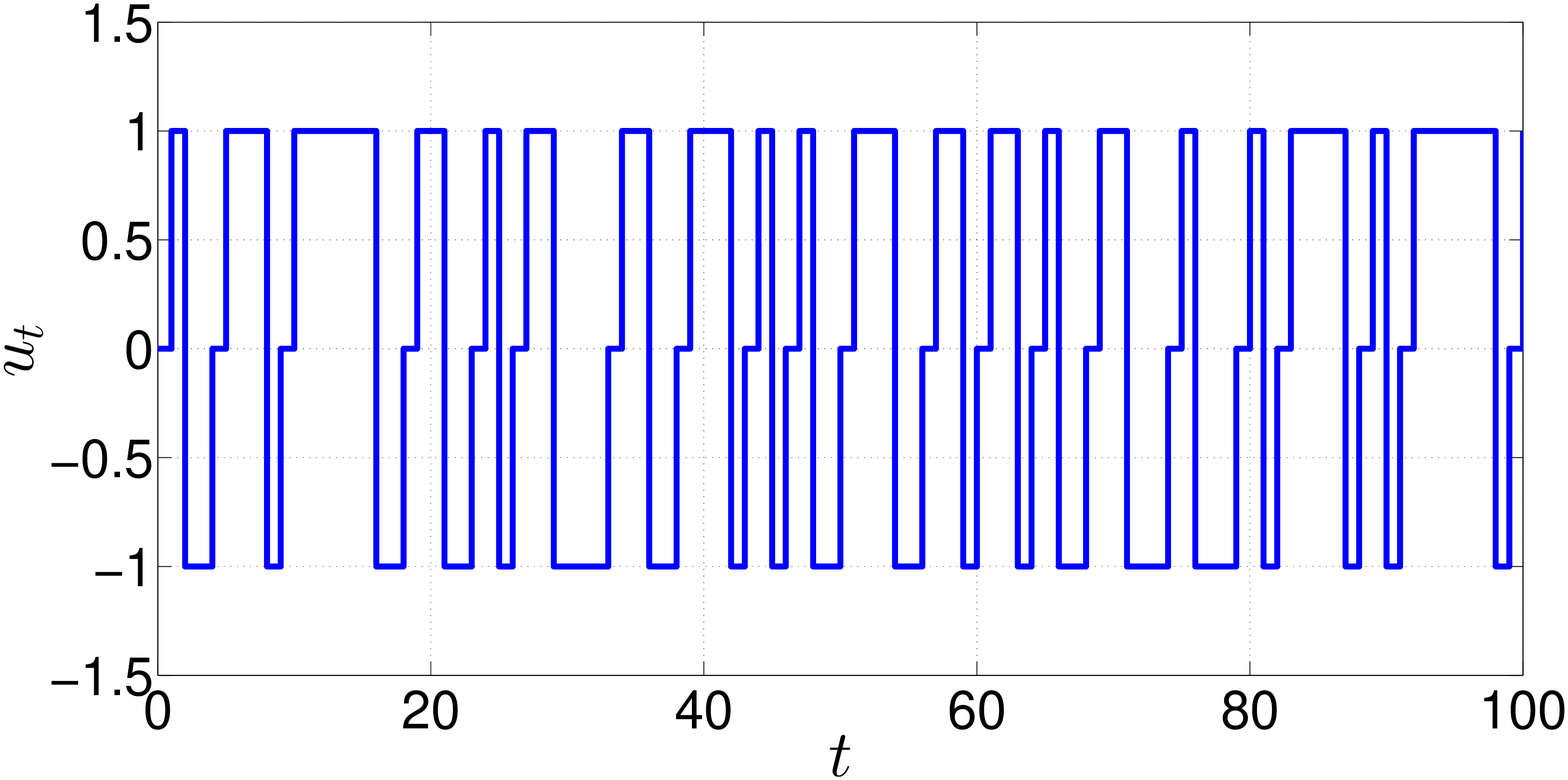}
\caption{Input realization, Example~\ref{ex: 2}.}}
\label{fig: input_det}
\end{figure}
%
%We generate $\{u_t ^j\}_{t=0} ^{N_{\op{sim}}}$ for each $v_j \in \Vcal_{\Pcal_\Ccal}$ ($N_{\op{sim}} = 10^2$) to compute the approximation \eqref{eq: var_score_app} for each $\Ical _F ^{(j)}$. Finally, the optimal input $u_{1:n_{\op{seq}}}$ is computed by running a Markov chain with $p^{\op{opt}} _u(u _{1:n_{m}})$ as stationary pmf, %\eqref{eq: prob14},
% where we discard the first $2 \cdot 10^6$ samples and keep the last $n_{\op{seq}} = 5 \cdot 10^3$ ones. A typical input realization is presented in Figure~\ref{fig: input_det}. In addition, we consider $K=100$, $M = 5\cdot 10^3$ and $N = 10^3$.
\begin{table}[t]
\caption{$h_m(\widehat{\Ical}_F )$, % for different input sequences,
 Example \ref{ex: 2}.}
\centering
\begin{tabular}{c|c}
  Input / $h_m(\widehat{\Ical}_F )$ & $\textsf{log } \textsf{det} (\widehat{\Ical}_F )$  \\%& $\textsf{tr} \left\{ (\widehat{\Ical}_F  )^{-1} \right\}$ \\
  \hline
  Optimal  & $\mathbf{25.34}(0.01)$ \\ %& $3.73 \cdot 10^{-5} (1.42 \cdot 10^{-7})$ \\
  %Optimal (tr)  & $21.56(0.04)$ & $2.68 \cdot 10^{-4} (5.85 \cdot 10^{-6})$ \\
  Binary        & $24.75(0.01)$ \\%& $6.70 \cdot 10^{-5} (2.54 \cdot 10^{-7})$ \\
  Uniform       & $24.38(0.01)$ %& $\mathbf{1.72 \cdot 10^{-4}} (1.42 \cdot 10^{-6})$ \\
\end{tabular}
\label{tab: inputs_2}
\end{table}
%As a benchmark, we use the same $n_{\op{seq}}$ input samples from the uniform and binary distribution presented in Example \ref{ex: 1}, which are employed to compute an approximation of $\Ical _F$ %the information matrix
% via \eqref{eq: var_score_app}.

Table \ref{tab: inputs_2} presents the results obtained for each input sequence, where \emph{Optimal} represents the result for the input sequence obtained from optimizing $\textsf{det} (\Ical_F ^{\op{app}} (\gamma))$. The $95\%$ confidence intervals are given as $\pm$ the value in the parentheses. From these data we conclude that the extended input design method outperforms the experiment results obtained for binary and uniformly distributed samples. Therefore, our new input design method can be successfully employed to design experiments for this nonlinear system. \fin

%Table \ref{tab: inputs_2} presents the results obtained for each input sequence, where \emph{Optimal} represents the result for the input sequence obtained from optimizing $\textsf{det} (\Ical_F ^{\op{app}} (\gamma))$. From these data we conclude that the extended input design method outperforms the experiment results obtained for binary and uniformly distributed samples. Therefore, our new input design method can be successfully employed to design experiments for the nonlinear system \eqref{eq:growth}. %Moreover, our result can be compared with those presented in \cite[Section 6]{gopaluni2011} when $h(\Ical_F ^{\op{app}}(\gamma)) = -\textsf{tr} \left\{ (\Ical_F ^{\op{app}} (\gamma) )^{-1} \right\}$.
%, for this particular example, the binary white noise seems to be the best input sequence. Indeed, the input design method proposed in this article tries to mimic the binary white noise excitation, which is clear from the numbers in Table \ref{tab: inputs}.
\end{exmp}
%%%%%%%%%%%%%%%%%%%%%%%%%%%%%%%%%%%%%%%%%%%%%%%%%%%%%%%%%%%%%%%%%%%%%%%%%%%%%%%%
\section{Conclusion}\label{sec: 5}
We have presented a new input design method for state space models, which extends existing input design approaches for nonlinear systems. The extension considers a more general model structure, and a new class for the input sequences. The method maximizes a scalar cost function of the information matrix, by optimizing the stationary pmf from which the input sequence is sampled. The elements in the feasible set of the stationary pmf are computed as a convex combination of its extreme points.

Under the assumption of a finite set of possible values for the input, we use graph theoretical tools to compute the information matrix as a convex combination of the information matrices associated with each extreme point. The information matrix for each extreme point is approximated using particle methods, where the information matrix is computed as the covariance of the score function. The numerical examples show that the extended input design method can be successfully used to design experiments for general nonlinear systems.

%We have presented a new input design method for state space models, which extends existing input design approaches for nonlinear systems. The extension considers a more general model structure, and a new class for the input sequences. The method maximizes a scalar cost function of the information matrix, by optimizing the stationary pmf from which the input sequence is sampled. The elements in the feasible set of pmf's are computed as a convex combination of its extreme points. Under the assumption of a finite set of possible values for the input, we use graph theoretical tools to compute the information matrix as a convex combination of the information matrices associated with each extreme point. The information matrix for each extreme point is approximated using particle methods, where the information matrix is computed as the covariance of the score function. The numerical examples show that the extended input design method can be successfully used to design experiments for general nonlinear systems.

In a future work we will combine the proposed input design technique with parameter estimation methods, which will allow to simultaneously estimate the parameters and the optimal input for a nonlinear SSM. We will also consider alternative methods based on Gaussian process models for information matrix estimation. This could improve the accuracy and the efficiency in the information matrix estimation method outlined in this paper.

Finally, as with most optimal input design methods, the one proposed in this contribution relies on knowledge of the true system. This difficulty can be overcome by implementing a robust experiment design scheme on top of it \citep{Rojas2007} or via an adaptive procedure, where the input signal is re-designed as more information is being collected from the system \citep{rojas2011adaptive}. This approach will be also addressed in a future work. %Future work on the subject will be focused on robust experiment design for nonlinear systems, where we will remove the assumption that the system parameters are known.

\begin{ack}                               % Place acknowledgements here.
The authors thank to Dr.\ Fredrik Lindsten for his comments to improve this article.
\end{ack}

\bibliographystyle{alpha}        % Include this if you use bibtex
%\bibliography{autosam}           % and a bib file to produce the
%\bibliography{autosam}
\bibliography{library,library2,dahlin}
                                 % bibliography (preferred). The
                                 % correct style is generated by
                                 % Elsevier at the time of printing.
%\appendix
%\section{A summary of Latin grammar}    % Each appendix must have a short title.
%\section{Some Latin vocabulary}         % Sections and subsections are supported
                                        % in the appendices.
\end{document}